\documentstyle[11pt,amstex,amssymb]{amsart}\textheight=19.8truecm\textwidth=14.1truecm
\begin{document}

\centerline{\Large\bf Concave functions of positive operators,}

\vskip 10pt \centerline{\Large\bf sums and congruences}

\vskip 15pt \centerline{{\sl \ Jean-Christophe Bourin \, and \,
Eun-Young Lee}} \vskip 5pt \centerline{{\small Daegu, January 28,
2007}}

\vskip 20pt \noindent {\small {\bf Abstract.}  Let $A$, $B$, $Z$ be
positive semidefinite matrices of same size and suppose $Z$ is
expansive, i.e.,  $Z\ge I$. Two remarkable inequalities are
$$
\Vert f(A+B)\Vert \le \Vert f(A)+f(B)\Vert \quad {\rm and} \quad
\Vert f(ZAZ)\Vert \le \Vert Zf(A)Z\Vert
$$
for all non-negative concave function $f$ on $[0,\infty)$  and all
symmetric norms $\|\cdot\|$ (in particular for all Schatten
$p$-norms). In this paper we survey several related results and we
show that these inequalities are two aspects of a unique theorem.
For the operator norm, our result also holds for operators on an
infinite dimensional Hilbert space.


Mathematical subjects classification:   15A60, 47A30}

\vskip 25pt\noindent {\large\bf 1. Introduction} \vskip 10pt A good
part of Matrix Analysis consists in establishing results for
Hermitian operators considered
 as generalized real numbers or functions.  In Section 2 we  recall two recent
 norm inequalities which are matrix versions of the obvious scalars
 inequalities
 $$
 f(za)\le zf(a) \quad {\rm and} \quad f(a+b) \le f(a)+f(b)
 $$
 for  non-negative concave functions $f$ on $[0,\infty)$  and
 scalars $a,\,b\ge0$ and $z\ge1$. In Section 3 we unify and generalize these
 norm inequalities. The norms considered are the symmetric (or unitarily invariant)
 norms. Such
 norms satisfy  $\Vert A\Vert = \Vert UAV\Vert$ for all $A$ and all
unitaries $U,\, V$. Here and in the sequel capital letters $A$,
$B,\dots,Z$ mean $n$-by-$n$ complex matrices, or operators on an
  $n$-dimensional Hilbert space. If $A$ is positive (semi-definite), resp.\ positive definite,
  we write $A\ge 0$, resp.\ $A>0$. If  $Z^*Z$ dominates the identity
  $I$, we say that $Z$ is expansive. As a corollary of Section 3 we
  have:

\vskip 10pt\noindent
 {\bf Theorem 1.1.}  {\it Let $\{A_i\}_{i=1}^m$ be positive and  let  $\{Z_i\}_{i=1}^m$ be expansive.
  Then, for all symmetric norms and all $p>1$,
$$
\left\| \sum Z_i^*A_i^pZ_i \right\| \le \left\|  \left(\sum
Z_i^*A_iZ_i \right)^p \right\|.
$$
} \noindent The sum is over the first $m$ integers. If $Z_i=I$ for
all $i$, this is a famous result of Ando and Zhan [1] and of Bhatia
and Kittaneh [3] in case of integer exponents. The very special case
$ {\mathrm Tr\,} (A_1^P+A_2^p) \le {\mathrm Tr\,} (A_1+A_2) ^p$ is
Mc-Carthy's inequality [13, p.\ 20].

 \vskip 20pt\noindent {\large\bf 2. Two norm inequalities}

 \vskip 10pt  The cone of positive operators
   is invariant under congruences $A\longrightarrow S^*AS$ and there
   are  several inequalities involving congruences  with a contraction $Z$ and  concave functions
   $f:[0,\infty)\rightarrow [0,\infty)$.  Brown-Kosaki's trace
   inequality states
   \begin{equation}
   {\mathrm Tr\,}f(Z^*AZ) \ge {\mathrm Tr\,}Z^*f(A)Z
   \end{equation}
   Actually a stronger result holds [8] (see also [4]), there exists a unitary $V$ such
   that
   \begin{equation}
   f(Z^*AZ) \ge VZ^*f(A)ZV^*.
   \end{equation}
   This means that the eigenvalues of $f(Z^*AZ)$ dominate those of
   $Z^*f(A)Z$.
    Further, if $f$ operator concave (equivalently operator monotone
   [10]), then  Hansen's operator inequality holds [9, 10],
   \begin{equation}
   f(Z^*AZ) \ge Z^*f(A)Z^*.
   \end{equation}
   What happens to (1), (2), (3) when $Z$ is no longer contractive but, in a opposite
   way, is  expansive ? It is
   obvious that (3) is reversed, meanwhile (2) can not be reversed
    though a non-trivial proof [4] shows that (1) is reversed:
     \begin{equation}
   {\mathrm Tr\,}f(Z^*AZ) \le {\mathrm Tr\,}Z^*f(A)Z.
   \end{equation}

   However a quite unexpected phenomena occurs: (1) can be extended to
   all Hermitians $A$ and all concave functions $f$ on the real line
   with  $f(0)\ge0$, but in (4) the assumption $A\ge0$ can not be
   dropped. The
good statement for expansive congruences requires positivity  and
involves
 symmetric (or unitarily invariant) norms $\Vert\cdot\Vert$.  We  have [5]:

 \vskip 10pt \noindent {\bf Theorem 2.1.} {\it Let
$f:[0,\infty)\longrightarrow [0,\infty)$ be a concave function. Let
$A\ge0$  and let $Z$ be expansive. Then, for all symmetric norms,
\begin{equation*}
\Vert f(Z^*AZ) \Vert \le \Vert Z^*f(A)Z\Vert.
\end{equation*}
}

 Besides these inequalities for congruences, there are nice
subadditivity results for concave functions
$f:[0,\infty)\longrightarrow [0,\infty)$  and sums of positive
operators. The most elementary one is a trace inequality companion
to (1) credited to Rotfel'd [12]: For $A,\, B\ge 0$,
\begin{equation}
   {\mathrm Tr\,}f(A+B) \le {\mathrm Tr\,} f(A) + f(B).
   \end{equation}
   The operator inequality $f(A+B) \le  f(A) + f(B)$ may not hold
   even if $f$ is operator concave. However, a remarkable subadditivity inequality
   related to (2) is shown in [2]: There exist unitaries $U,\, V$ such that
\begin{equation}
f(A+B) \le Uf(A)U^* + Vf(B)V^*.
\end{equation}
From this, it follows that the map
$$X\longrightarrow\Vert
f(|X|)\Vert$$
 is subadditive, see [7]. This was first noted by Uchiyama [14]. In case of the trace norm,
this extension of (5) is Rotfel'd's theorem [12]. Of course (6)
considerably strenghtens (5). Another improvement of (5), companion
to Theorem 1.1, is shown in [7]:

\vskip 10pt\noindent
 {\bf Theorem 2.2.}  {\it Let $A,\, B\ge 0$ and let  $f:[0,\infty)\longrightarrow[0,\infty)$ be a concave
function. Then, for all symmetric norms,
$$
\Vert f(A+B)\Vert \le \Vert f(A)+f(B)\Vert.
$$
}

\vskip -5pt \noindent  In case of the operator norm, Kosem [11] gave
a short proof. The general case is considerably more difficult: When
$f$ is operator concave, Theorem 2.2 had first been proved by Ando
and Zhan [1] by using integral representation of operator concave
functions and a delicate process. The proof given in [7] is much
more elementary.

\vskip 20pt\noindent {\large\bf 3. \ Combined result and proof}

\vskip 10pt We can naturally embodied Theorem 2.1 and Theorem 2.2
(and its version for sums of several operators) in a unique
statement:

\vskip 10pt\noindent
 {\bf Theorem 3.1.}  {\it Let $\{A_i\}_{i=1}^m$ be positive, let  $\{Z_i\}_{i=1}^m$ be expansive and let
 $f$ be a non-negative concave function on $[0,\infty)$. Then, for all symmetric norms,
$$
\left\|f\left(\sum Z_i^*A_iZ_i\right) \right\| \le \left\|  \sum
Z_i^*f(A_i)Z_i \right\|.
$$
}

\noindent Of course the sum is over the first $m$ integers. For
$m=1$ we get Theorem 2.1 and if $Z_i=I$ for all $i$ we get Theorem
2.2. There is no obvious way to derive Theorem 3.1 from Theorems
2.1, 2.2. However our proof is adapted from the proof of Theorem 2.1
which is itself partially based on the proof of Theorem 2.2.
Therefore our proof is rather elementary. We  use  Hansen's
inequality (3) for a quite elementary case equivalent to the fact
that $t\rightarrow 1/t$ is operator convex on $(0,\infty)$. We also
use some basic facts about symmetric norms. In particular we need
the following two facts for arbitrary $X, Y \ge0$.

\vskip 5pt\noindent 1) If $\| X\| \le \| Y\|$ for all symmetric
norms, then we also have $\| g(X)\| \le \| g(Y)\|$ for all symmetric
norms and all increasing convex functions
$g:[0,\infty)\longrightarrow[0,\infty)$.

\vskip 5pt\noindent 2) If $\| X\|_k \le \| Y\|_k$ for all Ky Fan
$k$-norms, then we also have $\| X\| \le \| Y\|$ for all symmetric
norms (Ky Fan's principle). Recall that the Ky Fan $k$-norms are the
sum of the $k$ largest singular values.

\vskip 5pt\noindent
 For this background we refer to any expository text such as [13].

\vskip 10pt\noindent
 {\bf Proof of Theorem 3.1.} The proof for an arbitrary $m$ is the
 same than the proof for $m=2$, so we consider the case of two expansive operators $X$,
 $Y$  and two positive operators $A$, $B$. The proof is divided in four steps.

\vskip 5pt Step 1.
  If $f$ is operator concave, the
 proof immediately follows from Theorem 2.1 and (3):
 \begin{align}
 \Vert f(X^*AX + Y^*BY) \Vert &\le   \Vert f(X^*AX) + f(Y^*BY) \Vert
 \notag
  \\
 &\le   \Vert X^*f(A)X + Y^*f(B)Y \Vert.
 \end{align}

 \vskip 10pt Step 2.
 Now consider a one to one convex function $g:[0,\infty)\longrightarrow
 [0,\infty)$ whose inverse function $f$ is operator concave. Since
 $g$  is onto there exist $A',\,B'\ge0$ such that $A=g(A')$ and
 $B=g(B')$; moreover $A'$ and $B'$ can be chosen arbitrarily so that  (7) can be read as
 $$
 \Vert f(X^*g(A)X + Y^*g(B)Y) \Vert \le   \Vert X^*AX + Y^*BY \Vert.
 $$
 Since $g$ is convex increasing we infer
\begin{equation}
 \Vert X^*g(A)X + Y^*g(B)Y \Vert \le   \Vert g(X^*AX + Y^*BY) \Vert.
 \end{equation}

\vskip 10pt
 Step 3.
 Now we extend (8) to the class of all non-negative convex functions
 on $[0,\infty)$ vanishing at $0$.  It
suffices to consider the Ky Fan $k$-norms $\Vert\cdot\Vert_{k}$.
Suppose that $g_1$ and $g_2$ both satisfy (8). Using the triangle
inequality and the fact that $g_1$ and $g_2$ are non-decreasing,
\begin{align*}
\Vert X^*(g_1+g_2)(A)X + &Y^*(g_1+g_2)(B)Y\Vert_{k} \\
&\le \Vert X^*g_1(A)X + Y^*g_1(B)Y\Vert_{k} + \Vert X^*g_2(A)X+Y^*g_2(B)Y\Vert_{k} \\
&\le \Vert g_1(X^*AX+Y^*BY)\Vert_{k} + \Vert g_2(X^*AX+Y^*BY)\Vert_{k} \\
&= \Vert (g_1+g_2)(X^*AX+Y^*BY)\Vert_{k},
\end{align*}
hence the set of functions satisfying to (8) is a cone. It is also
closed for point-wise convergence. Since any positive convex
function vanishing at $0$ can be approached by a positive
combination of angle functions at $a>0$,
$$
\gamma(t) =\frac{1}{2}\{ |t-a| + t-a\},
$$
it suffices to prove (8) for such a $\gamma$. By Step 2 it suffices
to approach $\gamma$ by functions whose inverses are operator
concave. We take (with $r>0$)
$$
h_r(t)=\frac{1}{2}\{ \sqrt{(t-a)^2 +r} + t -\sqrt{a^2 +r}\},
$$
whose inverse
$$
t-\frac{r/2}{2t+\sqrt{a^2 +r}-a} + \frac{\sqrt{a^2 +r}+a}{2}
$$
is operator concave since $1/t$ is operator convex on the positive
half-line. Clearly, as $r\to 0$, $ h_r(t)$ converges uniformly to
$\gamma$.

\vskip 10pt Step 4. Proof for any concave function
$f:[0,\infty)\longrightarrow[0,\infty)$. Again, it suffices to
consider the Ky Fan $k$-norms. This shows that we may and do assume
$f(0)=0$. Note that $f$ is necessarily non-decreasing. Hence, there
exists a rank $k$ spectral projection $E$ for $X^*AX+Y^*BY$,
corresponding to the $k$-largest eigenvalues
$\lambda_1(X^*AX+Y^*BY),\dots,\lambda_k(X^*AX+Y^*BY)$ of
$X^*AX+Y^*BY$, such that
$$
\Vert f(X^*AX+Y^*BY)\Vert_{k}=\sum_{j=1}^k
\lambda_j(f(X^*AX+Y^*BY))={\rm Tr\,} Ef(X^*AX+Y^*BY)E.
$$
Therefore, using a well-known property of Ky Fan norms, it suffices to show that
$$
{\rm Tr\,} Ef(X^*AX+Y^*BY)E \le {\rm Tr\,} E(X^*f(A)X+Y^*f(B)Y)E.
$$
This is the same as requiring that
\begin{equation}
{\rm Tr\,} E(X^*g(A)X+Y^*g(B)Y)E \le {\rm Tr\,} Eg(X^*AX+Y^*BY)E
\end{equation}
for all non-positive convex functions $g$ on $[0,\infty)$ with
$g(0)=0$. Any such function can be approached by a combination of
the type
\begin{equation*}
g(t)=\lambda t + h(t)
\end{equation*}
for a scalar $\lambda <0$ and some non-negative convex function $h$
vanishing at 0. Hence, it suffices to show that (9) holds
 for $h(t)$.
We have
\begin{align*}
{\rm Tr\,} E(X^*h(A)X+Y^*h(B)Y)E &= \sum_{j=1}^k \lambda_j(E(X^*h(A)X+Y^*h(B)Y)E) \\
&\le \sum_{j=1}^k \lambda_j(X^*h(A)X+Y^*h(B)Y) \\
&\le \sum_{j=1}^k \lambda_j(h(X^*AX+Y^*BY)) \quad {\rm (by\ Step\ 3)} \\
&= \sum_{j=1}^k \lambda_j(Eh(X^*AX+Y^*BY)E) \\
&= {\rm Tr\,} Eh(X^*AX+Y^*BY)E
\end{align*}
where the second equality follows from the fact that $h$ is
non-decreasing and hence $E$ is also a spectral projection of
$h(X^*AX+Y^*BY)$ corresponding to the $k$ largest eigenvalues.
 \qquad $\Box$

\vskip 10pt\noindent
 {\bf Corollary 3.2.}  {\it Let $g:[0,\infty)\longrightarrow[0,\infty)$ be a convex function with
  $g(0)=0$. Let $\{A_i\}_{i=1}^m$ be positive and let  $\{Z_i\}_{i=1}^m$ be expansive. Then, for all symmetric norms,
$$
\left\| \sum Z_i^*g(A_i)Z_i \right\| \le \left\|  g\left(\sum
Z_i^*A_iZ_i \right) \right\|.
$$
}

\noindent This corollary is proved in step 3. It can also be derived
from Theorem 3.1 by using the first fact recalled before the proof.
When $Z_i=I$ for all $i$ this is a remarkable result of Kosem [11].
Note that Theorem 1.1 is a special case of Corollary 3.2.

\vskip 10pt There are several well-known results involving sums of
congruences. But, in contrast with Theorem 3.1 and its corollary,
these results deal with strong contractive assumptions. We give an
example generalizing (2):

\vskip 10pt\noindent
 {\it Let $\{A_i\}_{i=1}^m$ be positive and
$\{Z_i\}_{i=1}^m$ such that $\sum Z^*_iZ_i \le I$. If
 $f$ is a monotone concave function on $[0,\infty)$, $f(0)\ge0$, then,
$$
f\left(\sum Z_i^*A_iZ_i\right)  \ge   V\left(\sum Z_i^*f(A_i)Z_i
\right)V^*.
$$
 for some unitary $V$.}

\vskip 10pt\noindent Such an inequality is connected to  Jensen type
inequalities for compressions or positive unital linear maps, see
[4], [2].

\vskip 10pt We conclude with a discussion of the extension of our
results for operators on  an infinite dimensional Hilbert space
$\mathcal{H}$. In the proof of Theorem 3.1, (7) is derived from a
version of Hansen's inequality (3) involving congruences with
expansive operators $X, Y$. This Hansen's inequality remains valid
in the infinite dimensional setting if one requires that $X, Y$ are
both expansive and invertible. Consequently we have the following
result for the usual operator norm $\|\cdot\|_{\infty}$.

\vskip 10pt\noindent
 {\it Theorem 3.1 and its corollaries are still valid for $\|\cdot\|_{\infty}$
 and operators on $\mathcal{H}$ when $Z$ and
$\{Z_i\}_{i=1}^m$ are expansive and invertible.}

\vskip 10pt\noindent This statement is meaningful. The original
proof ([5]) of Theorem 2.1 was unsuccessful to cover the infinite
dimensional setting. Concerning inequality (3), the original
statement is in the framework of the spectral order in a semi-finite
von Neumann algebra. It is also possible to give a version for
operators on $\mathcal{H}$ by adding a $rI$ term in the RHS with
$r>0$ arbitrarily small, see [6] pp.\ 11-15. By arguing as in [2] we
then obtain:

\vskip 10pt\noindent
 {\it Let $A, B \ge 0$ on $\mathcal{H}$ and let $r>0$. If
 $f$ is a monotone concave function on $[0,\infty)$, $f(0)\ge0$, then,
$$
f(A+B)  \le  Uf(A)U^* + Vf(B)V^* +rI
$$
 for some unitaries $U,\,V$.}

\newpage
 \vskip 10pt {\bf References}

\noindent {\small \vskip 5pt\noindent [1] T.\ Ando and X.\ Zhan,
Norm inequalities related to operator monotone functions, Math.\
Ann.\ 315 (1999) 771-780.

\vskip 5pt\noindent [2] J.\ S.\ Aujla and J.-C.\ Bourin, Eigenvalue
inequalities for convex and log-convex functions, Linear Alg.\
Appl., in press (2007).

 \vskip 5pt\noindent
[3] R.\ Bhatia and F.\ Kittaneh, Norm inequalities for positive
operators, Lett.\ Math.\ Phys.\ 43 (1998) 225-231.

 \vskip 5pt\noindent [4] J.-C.\ Bourin, Convexity or  concavity inequalities for
 Hermitian operators, Math.\ Ineq.\ Appl., 7 (4) (2004) 607-620.

 \vskip 5pt\noindent [5] J.-C.\ Bourin, A concavity inequality for symmetric
norms, Linear Alg.\ Appl., 413 (2006) 212-217.

\vskip 5pt\noindent [6] J.-C.\ Bourin, Compressions, Dilations and
Matrix Inequalities, RGMIA monograph, Victoria University,
Melbourne, 2004, Available
from:$<$http://rgmia.vu.edu.au/monographs$>$

 \vskip 5pt\noindent [7] J.-C.\ Bourin and M.\ Uchiyama, A matrix subadditivity
 inequality for  $f(A+B)$ and $f(A)+f(B)$, Linear Alg.\ Appl., in press.

\vskip 5pt\noindent [8] L.\ G.\ Brown and H. Kosaki, Jensen's
inequality in semi-finite von Neuman algebras, J. Operator theory 23
(1990) 3-19.

\vskip 5pt\noindent [9] F.\ Hansen, An operator inequality, Math.\
Ann.\ 246 (1980) 249-250.

\vskip 5pt\noindent [10] F.\ Hansen and G.\ K.\ Pedersen, Jensen's
inequality for operators and Lowner's theorem, Math.\ Ann.\ 258
(1982) 229-241.

 \vskip 5pt\noindent [11] T.\ Kosem,
Inequalities between $\Vert f(A+B)\Vert$ and $\Vert f(A)+f(B)\Vert$,
Linear Alg.\ Appl., 418 (2006) 153-160.

\vskip 5pt\noindent [12] S. Ju. Rotfel'd, The singular values of a
sum of completely continuous operators, Topics in Mathematical
Physics, Consultants Bureau, Vol.\ 3, 1969, 73-78.

\vskip 5pt\noindent [13] B.\ Simon, Trace Ideals and Their
Applications, LMS lecture note 35, Cambridge Univ.\ Press,
Cambridge, 1979.

 \vskip 5pt\noindent [14] M.\
Uchiyama, Subadditivity of eigenvalue sums, Proc.\ Amer.\ Math.\
Soc., 134 (2006) 1405-1412.
 }

\vskip 10pt \centerline{Jean-Christophe Bourin }
 \centerline{E-mail: bourinjc@@club-internet.fr }
  \centerline{Department of mathematics}
  \centerline{Kyungpook National University}
  \centerline{ Daegu 702-701, Korea}

  \vskip 10pt \centerline{ Eun-Young Lee }
 \centerline{E-mail: eylee89@@knu.ac.kr}
  \centerline{Department of mathematics}
  \centerline{Kyungpook National University}
  \centerline{ Daegu 702-701, Korea}

\end{document}